\documentclass[12pt,reqno]{amsart}

\usepackage[utf8]{inputenc}
\usepackage[english]{babel}
\usepackage{amsmath,amssymb,amsfonts,amsbsy,amsthm,amscd,latexsym,graphicx}
\usepackage{empheq,textcomp}
\usepackage{framed}

\theoremstyle{definition}
\newtheorem{theorem}{Theorem} [section]
\newtheorem{corollary}[theorem]{Corollary}
\newtheorem{lemma}[theorem]{Lemma}
\newtheorem{proposition}[theorem]{Proposition}
\newtheorem{definition}[theorem]{Definition}

\newtheorem{remark}[theorem]{Remark}
\newtheorem{example}[theorem]{Example}
\newtheorem{conjecture}[theorem]{Conjecture}
\numberwithin{equation}{section}

\newcommand{\C}{\mathbb{C}}
\newcommand{\Cc}{{\mathcal{C}}}
\newcommand{\Dc}{{\mathcal{D}}}
\newcommand{\Ec}{{\mathcal{E}}}
\newcommand{\Fc}{{\mathcal{F}}}
\newcommand{\Gc}{{\mathcal{G}}}
\newcommand{\N}{\mathbb{N}}
\newcommand{\R}{\mathbb{R}}
\newcommand{\Rc}{{\mathcal{R}}}
\newcommand{\T}{\mathbb{T}}
\newcommand{\Yc}{{\mathcal{Y}}}
\newcommand{\Z}{\mathbb{Z}}

\newcommand{\eps}{\varepsilon} %Better epsilon symbol
\newcommand{\tEc}{{\widetilde{\mathcal{E}}}}
\newcommand{\tFc}{{\widetilde{\mathcal{F}}}}
\newcommand{\tRc}{{\widetilde{\mathcal{R}}}}
\newcommand{\te}{{\widetilde{e}}}
\newcommand{\tr}{{\widetilde{r}}}
\newcommand{\tx}{{\widetilde{x}}}

\newcommand{\Eq}{\, = \,}
\newcommand{\Ne}{\, \ne \,}
\newcommand{\Le}{\, \le \,}
\newcommand{\Lt}{\, < \,}
\newcommand{\Plus}{\, + \,}
\newcommand{\qeddef}{{\quad $\diamondsuit$}}
\newcommand{\qeddeff}{{\qquad \diamondsuit}}

\newcommand{\bigabs}[1]{\bigl|#1\bigr|}
\newcommand{\biggabs}[1]{\biggl|#1\biggr|}
\newcommand{\ip}[2]{\langle#1,#2\rangle}
\newcommand{\biggip}[2]{\biggl\langle #1, \, #2 \biggr\rangle}
\newcommand{\norm}[1]{\|#1\|}
\newcommand{\biggnorm}[1]{\biggl\|#1\biggr\|}

\newcommand{\Bigparen}[1]{\Bigl(#1\Bigr)}
\newcommand{\biggparen}[1]{\biggl(#1\biggr)}
\newcommand{\range}{{\textup{range}}}
\newcommand{\set}[1]{\{#1\}}
\newcommand{\bigset}[1]{\bigl\{#1\bigr\}}
\newcommand{\Bigset}[1]{\Bigl\{#1\Bigr\}}
\newcommand{\biggset}[1]{\biggl\{#1\biggr\}}
\newcommand{\clspan}{{\overline{\text{span}}}}

\newcommand{\el}{{\ell^1}}

\begin{document}

\title{$\ell^1$-Bounded Sets}

\author{Christopher Heil and Pu-Ting Yu}

\address{School of Mathematics, Georgia Institute of Technology,
Atlanta, GA 30332, USA}

\email{heil@math.gatech.edu}
\email{pyu73@gatech.edu}

\thanks{Acknowledgement:
This research was partially supported by a grant from the
Simons Foundation.}

\date{April 30, 2023}

\begin{abstract}
A subset $M$ of a separable Hilbert space $H$ is
\emph{$\el$-bounded} if there exists a Riesz basis
$\Fc = \set{e_n}_{n \in \N}$ for $H$ such that
$\sup_{x \in M} \sum_{n \in \N} |\ip{x}{e_n}| < \infty.$
A similar definition for \emph{$\ell^1$-frame-bounded sets}
is made by replacing Riesz bases with frames.
This paper derives properties of $\ell^1$-bounded sets,
operations on the collection of $\el$-bounded sets,
and the relation between $\el$-boundedness and $\el$-frame-boundedness.
Some open problems are stated, several of which have intriguing implications.
\end{abstract}

\maketitle

\section{Introduction}
A sequence $\set{x_n}_{n \in \N}$ in a separable Hilbert space $H$ is a
\emph{Riesz basis} for $H$ if it is the image of an orthonormal basis
under a bounded invertible linear operator.
We refer to \cite{Chr16}, \cite{Hei11}, and \cite{You01}
for relatively recent texts that discuss Riesz bases
and the related topic of \emph{frames}.

The Wiener algebra is a classical topic in harmonic analysis.
The issue there is to obtain conditions that imply summability of the
Fourier coefficients of a function in $L^1(\T).$
We will consider certain generalizations of the Wiener algebra in
separable Hilbert spaces.
Instead of focusing on the exponential family $\set{e^{2\pi inx}}_{n \in \Z},$
we work with a Riesz basis $\Ec = \set{e_n}_{n \in \N}$
for a separable Hilbert space $H.$
We define the $\ell^1$-norm of a subset $M$ of $H$ with respect to $\Ec$ to be
$\sup_{x \in M} \sum_{n \in \N} \, |\ip{x}{e_n}|.$
Following \cite{HH13},
we say that $M$ is $\el$-bounded if its $\el$-norm
with respect to \emph{some} Riesz basis is finite.

Frames, which are generalizations of Riesz bases,
were first introduced by Duffin and Schaeffer \cite{DS52}
in their study of non-harmonic Fourier series.
A sequence $\set{x_n}_{n \in \N}$ in a separable Hilbert space $H$
is a \emph{frame} if there exist positive constants $A \le B,$
called \emph{frame bounds}, such that 
\begin{equation} \label{framebound_eq}
A \, \norm{x}^2
\Le \sum_{n=1}^\infty |\ip{x}{x_n}|^2
\Le B \, \norm{x}^2,
\qquad\text{for all } x \in H.
\end{equation}
If we can take $A=B=1$ then we say that $\set{x_n}_{n \in \N}$
is a \emph{Parseval frame}.
The significance of equation \eqref{framebound_eq} is that it endows frames
with basis-like properties.
In particular, there exists at least one sequence $\set{y_n}_{n \in \N}$
in $H$ such that every $x$ in $H$ can be expressed as
\begin{equation} \label{frameseries_eq}
x \Eq \sum_{n=1}^\infty \, \ip{x}{y_n} \, x_n,
\end{equation}
where the sum converges in the norm of $H.$
A sequence $\set{x_n}_{n \in \N}$ is a Riesz basis if and only if it is a
nonredundant frame (which means that the removal of any element from
the sequence leaves an incomplete set).
The idea of $\ell^1$-boundedness can be formulated for frames as
well as Riesz bases.

In this paper, we study $\el$-boundedness, including topological features,
operations on the collection of $\el$-bounded sets, and the relation
between $\el$-boundedness and $\el$-frame-boundedness.

We begin in Section \ref{prelim_sec} with a presentation of
notation and terminology and a review of some basic properties
and characterizations of Riesz bases and frames.

In Section~\ref{H_space}, we study the space $H_\Ec$ consisting of
all elements of $H$ whose coefficients are
absolutely summable with respect to a particular Riesz basis $\Ec.$
We then consider $\el$-bounded sets in detail in
Section~\ref{ell1_bdd_sets_sec}.
Intriguingly, the question of whether the collection of $\el$-bounded
sets is closed under finite unions seems to be a very difficult question.
Among other results, we prove that this collection is closed under
finite unions if and only if it is closed under finite sums,
and that if such closure holds then several remarkable implications
follow, including that $\el$-boundedness would in that case be
equivalent to $\el$-frame boundedness.
We therefore conjecture that the collection of $\el$-bounded sets
is not closed under finite sums or unions.

Finally, in Section \ref{absolutelyconvergent_sec}
we use our results to study $p$-convergent frames,
whose frame series are $p$-summable.
In particular we prove that $p$-convergent frames exist for $1 < p \le 2,$
but there are no absolutely convergent frames.

We focus exclusively in this paper on complex Hilbert spaces.
Some of our results carry over to real Hilbert spaces,
but the proofs of several results ultimately rely on the
Polar Decomposition of bounded linear operators on complex Hilbert spaces.
We leave the question of whether those results have analogues
in the real setting as another open question.

\subsection*{Acknowledgements}
We thank Bernhard Haak and Markus Haase for formulating the definition of
$\ell^1$-bounded sets, and for bringing the question of
whether $\ell^1$-bounded sets are closed under finite unions
to our attention.

\section{Preliminaries} \label{prelim_sec}
Throughout this paper, $H$ denotes a separable infinite-dimensional
complex Hilbert space equipped with inner product $\ip{\cdot}{\cdot}$
and induced norm $\norm{\cdot}.$
Unless stated otherwise, $\set{x_n}$ will denote a countable sequence
indexed by the natural numbers $\N.$
When needed for clarity, we will write $\set{x_n}_{n \in \N}.$
We let $\ell^2$ denote the usual space of square-summable sequences of scalars.

If $S$ is a closed subspace of $H,$ then $P_S$ will denote the
orthogonal projection of $H$ onto $S.$
We refer to a bounded linear invertible operator that maps one
Hilbert space onto another as a \emph{topological isomorphism}.
We say that $H$ \emph{embeds} into a Hilbert space $K$ if
$H \subseteq K$ and $\ip{x}{y}_K = \ip{x}{y}_H$ for all $x,$ $y \in H.$

\smallskip
\subsection{$\ell^1$-Bounded Sets}

A sequence $\set{x_n}$ in $H$ is a
\emph{Bessel sequence} if it satisfies at least the upper inequality
in the definition of a frame.
That is, there exists some constant $B > 0$ such that
$\sum |\ip{x}{x_n}|^2 \le B \, \norm{x}^2$
for all $x \in H.$
We call $B$ an \emph{upper frame bound} for $\set{x_n}.$
Similarly, we say that a sequence possesses a lower frame bound
if there exists a constant $A>0$ such that
$\sum |\ip{x}{x_n}|^2 \ge A \, \norm{x}^2$ for all $x \in H.$

\begin{definition} \label{ell1bounded_def}
(a) If $\Ec = \set{e_n}$ is a Riesz basis or frame for $H$
and $x \in H,$ then the \emph{$\ell^1$-norm of $x$ with respect to $\Ec$} is
$$\norm{x}_{1,\Ec} \Eq \sum_n |\ip{x}{e_n}|.$$

(b) A subset $M$ of $H$ is \emph{$\ell^1$-bounded}
if there exists a Riesz basis $\Ec = \set{e_n}$ for $H$ such that
\begin{equation*} %\label{ell1boundeddef_eq}
\sup_{x \in M} \norm{x}_{1,\Ec} \Lt \infty.
\end{equation*}
Similarly, $M$ is \emph{$\ell^1$-frame-bounded}
if there exists a frame $\Fc = \set{x_n}$ for $H$ such that 
$$\sup_{x \in M} \norm{x}_{1,\Fc} \Lt \infty.$$

(c) Given a Riesz basis or frame $\Ec$ for $H,$ we let 
$H_\Ec$ denote the collection of all $x \in H$ such that
$\norm{x}_{1,\Ec}$ is finite.
That is,
$$H_\Ec \Eq \bigset{x \in H : \norm{x}_{1,\Ec} < \infty}. \qeddeff$$
\end{definition}

Although the notions in Definition \ref{ell1bounded_def} seem natural,
we are not aware of much literature on the subject, at least as is
related to the directions taken in this paper.
In particular, the space $H_\Ec$ with $\Ec$ being an orthonormal basis
was studied in \cite{HL08}
because of certain connections between operator theory
and modulation spaces in time-frequency analysis.
Some $\ell^p$ analogues of $H_\Ec$ were considered
by Han, Li, and Tang in \cite{HLT11}.
The idea of $\ell^1$-bounded and $\ell^1$-frame bounded sets
was introduced and studied by Haak and Haase in \cite{HH13}.
In addition, they obtained results on
\emph{$\ell^1$-bounded operators},
while we will focus only on $\ell^1$-bounded sets.

\smallskip
\subsection{Some Properties of Riesz Bases and Frames}

If $\Fc = \set{x_n}$ is a frame for $H,$ then its \emph{frame operator}
$$Sx \Eq \sum_n \ip{x}{x_n} \, x_n, \qquad\text{for } x \in H,$$
is a positive bounded linear invertible mapping of $H$ onto itself.
The \emph{canonical dual frame} of $\Fc$ is
$\tFc = \set{\tx_n}$ where $\tx_n = S^{-1}x_n$ for every $n.$
This is a frame for $H,$ and for every $x \in H$ we have the
\emph{reconstruction formulas}
\begin{equation} \label{reconstruction_eq}
x \Eq \sum_n \ip{x}{\tx_n} \, x_n \Eq \sum_n \ip{x}{x_n} \, \tx_n.
\end{equation}
The series in equation \eqref{reconstruction_eq}
converge unconditionally (that is, regardless of the ordering of the
index set) for every $x \in H.$
If $\Fc$ is a Parseval frame then $S = I,$ the identity operator.
If $\Fc$ is any frame, then $S^{-1/2}(\Fc)$ is a Parseval frame.

Any sequence $\set{y_n}$ in $H$ such that
$x = \sum \ip{x}{y_n} \, x_n$ for $x \in H$
is called an \emph{alternative dual} for $\Fc.$
If an alternative dual is also a frame, then it is called an
\emph{alternative dual frame}.

If $\Fc$ is a Riesz basis for $H,$ then the canonical dual frame
is the unique alternative dual for $\Fc.$
Further, in this case $\tFc$ is itself a Riesz basis, called the
\emph{dual basis} for $\Fc.$
This dual basis is the unique sequence that is \emph{biorthogonal}
to $\Fc$ in the sense that $\ip{x_m}{\tx_n} = \delta_{mn}.$

\smallskip
\subsection{Characterizations of Frames}

One proof of the following result can be found in \cite[Cor.~8.30]{Hei11}.

\begin{theorem} \label{frame_image_Riesz}
A sequence $\set{x_n}$ is a frame for $H$ if and only if there exists a
Riesz basis $\set{e_n}$ for $H$ and a bounded linear surjective operator
$T \colon H \to H$ such that $T(e_n) = x_n$ for every $n.$
\qeddef\end{theorem}

The next theorem is a version of the \emph{Naimark duality principle},
and states that every frame for $H$ is the orthogonal projection
of some Riesz basis for some larger Hilbert space.
The theorem appears to have been first stated explicitly in this form
by Han and Larson \cite{HL00}.

\begin{theorem} \label{Nai_Du}
A sequence $\set{x_n}$ is a frame for $H$ if and only if there exists
a Hilbert space $K$ such that $H$ embeds into $K$
and there is a Riesz basis $\set{e_n}$ for $K$
such that $x_n = P_H e_n$ for $n \in \N.$

Furthermore, by identifying $H$ with $H \times \set{0}$ we can construct
$K$ to have the form $K = H \times N^\perp,$
where $N$ is a closed subspace of $\ell^2$ that
is topologically isomorphic to $H.$
\qeddef\end{theorem}

\smallskip
\subsection{Operator Decompositions and Frames}

Casazza states in \cite{Cas98} that a well-known result in operator theory
is that any bounded linear operator $T \colon H \to H$ can be
written as $T = a(U_1 + U_2 + U_3)$ where $a$ is a nonnegative scalar
and $U_1,$ $U_2,$ and $U_3$ are unitary, and he provides a proof of
this result.
A corollary is that every frame $\set{x_n}$ can be written in the
form $x_n = a(e_n + f_n + g_n)$ for some orthonormal bases
$\set{e_n},$ $\set{f_n},$ and $\set{g_n}$ for $H.$
These facts hold for complex Hilbert spaces;
there are analogous results for real Hilbert spaces but more
unitary operators or orthonormal bases are required in that setting.

We will require a variation on Casazza's construction,
writing a bounded linear operator on a complex Hilbert space
as a multiple of a sum of a unitary operator and a topological isomorphism.
We include the proofs for completeness, but note that our
Lemma \ref{topisodecomp_lemma} is similar to \cite[Prop.~1]{Cas98},
while the proof of our Theorem \ref{Decom_operator}
is essentially the proof of \cite[Prop.~7]{Cas98}.

\begin{lemma} \label{topisodecomp_lemma}
If $T \colon H \to H$ is a topological isomorphism with $\norm{T} \le 1,$
then there exist unitary operators $U$ and $V$ such that
$T = \frac12 (U + V).$
\end{lemma}
\begin{proof}
Using the Polar Decomposition of bounded linear operators on
complex Hilbert spaces (see \cite[Thm.~VIII.3.11]{Con90}),
we can write $T = WP$ where $W$ is a unitary operator and $P = (T^*T)^{1/2}$
(we have that $W$ is unitary here because $T$ is a topological isomorphism).
Since $\norm{P} \le 1,$ we can define $Q = (I - P^2)^{1/2}.$
If we let $V = P + iQ,$ then
$P = \frac12 (V + V^*)$ and $T = \frac12 (WV + WV^*).$
Since $Q$ commutes with all operators that commute with $I - P^2,$
which includes $P$ itself, we have that $V^*V = I = VV^*.$
Therefore $V$ is unitary.
\end{proof}

\begin{theorem} \label{Decom_operator}
If $T \colon H \to H$ is a bounded linear operator,
then there exist a scalar $a \ge 0,$
a unitary operator $U,$
and a topological isomorphism $S$
such that $T = a(U + S).$
\end{theorem}
\begin{proof}
Fix $0 < \eps < 1,$ and let
$$W \Eq \frac34 I \Plus \biggparen{\frac{1-\eps}4} \, \frac{T}{\norm{T}}.$$
Then $\norm{I-W} < 1,$ so $W$ is a topological isomorphism.
By Lemma \ref{topisodecomp_lemma}, there exist unitary operators
$U$ and $V$ such that $W = \frac12 (U + V).$
Consequently,
$$T
\Eq \frac{4 \, \norm{T}}{1 - \eps} \Bigparen{W - \frac34 I}
\Eq \frac{2 \, \norm{T}}{1 - \eps}
    \biggparen{U \Plus \Bigparen{V - \frac32 I}}
\Eq a (U + S)$$
where $a = \frac{2 \, \norm{T}}{1 - \eps}$ and $S = V - \frac32I.$
Finally, since $\norm{I - (-\frac12 S)} < 1,$
we know that $-\frac12 S,$ and therefore $S$ itself,
is a topological isomorphism.
\end{proof}

It is shown in \cite{Cas98}
that every frame is a sum of finitely many Riesz bases.
We observe next that this decomposition extends to Bessel sequences.
Another proof of Theorem \ref{Bes_decomposition}
was given by Dehghan and Mesbah \cite{DM15}, based on a deep result in
$C^*$-algebras by Kadison and Pederson \cite[Thm.~1]{KP85}.

\begin{theorem} \label{Bes_decomposition}
If $\Fc = \set{x_n}$ is a Bessel sequence in $H,$ then there exist
Riesz bases $\set{y_n}$ and $\set{z_n}$ for $H$ such that
$x_n = y_n + z_n$ for every $n.$
\end{theorem}
\begin{proof}
Let $\set{x_n}$ be a Bessel sequence and let $\set{e_n}$
be an orthonormal basis for $H.$
Because $\set{x_n}$ is Bessel, the mapping
$$Tx \Eq \sum_{n=1}^\infty \ip{x}{e_n} \, x_n,
\qquad\text{for } x \in H,$$
is bounded and linear, and $T$ sends $e_n$ to $x_n$ for every $n.$
Let $T = a(U + S)$ be the decomposition given
by Theorem \ref{Decom_operator},
and let $y_n = a Ue_n$ and $z_n = a Se_n.$
Then $\set{y_n}$ and $\set{z_n}$ are each Riesz bases for $H$
(in fact, $\set{y_n}$ is a multiple of an orthonormal basis),
and we have $x_n = y_n + z_n$ for every $n.$
\end{proof}

\smallskip
\subsection{A Frame With $H_\Fc = \set{0}$}

We will see in Section \ref{H_space} that $H_\Ec$ is a dense subspace
of $H$ when $\Ec$ is a Riesz basis.
However, Han, Li, and Tang showed in \cite{HLT11} that $H_\Fc$
can be as small as $\set{0}$ when $\Fc$ is a frame.
Their construction is simple, so we include it for completeness.

\begin{example} \label{empty_HF_space}
Let $\set{e_n}$ be an orthonormal basis for $H.$
For each $n,$ let $(a_{jn})_{j \in \N}$
be a sequence of scalars such that $\sum_j |a_{jn}| = \infty$
and $\sum_j |a_{jn}|^2 = 1.$
Let $\Fc = \set{a_{jn} e_n}_{j, n \in \N}.$
If $x \in H,$ then
$$\sum_{j,n} |\ip{x}{a_{jn}e_n}|^2
\Eq \sum_n |\ip{x}{e_n}|^2 \, \sum_j |a_{jn}|^2
\Eq \sum_n |\ip{x}{e_n}|^2
\Eq \norm{x}^2.$$
Therefore $\Fc$ is a Parseval frame.
On the other hand, if $x \ne 0$ then
$$\sum_{j,n} |\ip{x}{a_{jn}e_n}|
\Eq \sum_n |\ip{x}{e_n}| \, \sum_j |a_{jn}|
\Eq \sum_n |\ip{x}{e_n}| \cdot \infty
\Eq \infty.$$
Hence $H_\Fc = \set{0}.$
\qeddef\end{example}

\smallskip
\subsection{Convergence of Dual Sequences}

It is known that, in at least some cases, convergence of the series
$\sum \, \ip{x}{y_n} \, x_n$ implies the convergence of
$\sum \, \ip{x}{x_n} \, y_n.$
Combining \cite[Lem.~3.1]{SB13} with \cite[Cor.~3.10]{SB11}
gives us the following result.

\begin{lemma} \label{equiv_conv_alt_syn_dual}
Let $\set{x_n}$ and $\set{y_n}$ be two sequences in $H.$
\begin{enumerate}
\item[\textup{(a)}]
$\sum \, \ip{x}{y_n} \, x_n$ converges unconditionally for every $x \in H$
if and only if $\sum \, \ip{x}{x_n} \, y_n$ converges unconditionally
for every $x \in H.$

\medskip
\item[\textup{(b)}]
If the equivalent conditions of statement (a) hold, then
$x = \sum \, \ip{x}{y_n} \, x_n$ for every $x \in H$
if and only if $x = \sum \, \ip{x}{x_n} \, y_n$
for every $x \in H.$ \qeddef
\end{enumerate}
\end{lemma}

\section{$H_\Ec$ Spaces} \label{H_space}

\subsection{Examples}

\begin{example} \label{l1_example}
If $H = \ell^2$ and $\Ec$ is the standard basis for $\ell^2,$
then $H_\Ec = \ell^1.$
\qeddef\end{example}

\begin{example}
If $H = L^2(\T)$ then the trigonometric system
$\Ec = \set{e^{2\pi inx}}_{n \in \Z}$ is an orthonormal basis for $L^2(\T).$
In this case $H_\Ec = A(\T),$ the Wiener algebra consisting of
those functions in $L^2(\T)$ with summable Fourier coefficients.
\qeddef\end{example}

\begin{example}
Let $H = L^2(\R),$ and let $\phi(x) = e^{-x^2}$ be the Gaussian function.
Fix $a,$ $b > 0$ with $ab < 1,$
and for $m,$ $n \in \Z$ let
$\phi_{nk}(x) = e^{2\pi ibnx} \phi(x-ka)$
be a time-frequency shift of $\phi.$
In this case the \emph{Gabor system} $\Gc = \set{\phi_{nk}}_{n,k \in \Z}$
is a frame for $L^2(\R).$
The space
$$H_\Gc
\Eq \Bigset{f \in L^2(\R) : \sum_{n,k \in \Z} |\ip{f}{\phi_{nk}}| < \infty},$$
consisting of functions whose Gabor coefficients are summable,
equals the \emph{Feichtinger algebra}, which is usually denoted by
$M^1(\R)$ or $S_0(\R).$
The Feichtinger algebra plays a central role in time-frequency analysis;
we refer to \cite{Gro01} for details.
\qeddef\end{example}

\smallskip
\subsection{Topological Features}

We will consider some topological properties of the space $H_\Ec.$

First, we summarize some basic facts about $H_\Ec$ for the case
where $\Ec$ is a Riesz basis.
  
\begin{proposition} \label{iso_ell1}
Assume $\Ec = \set{e_n}$ is a Riesz basis for $H,$
and let $\set{\te_n}$ be its dual basis.
\begin{enumerate}
\item[\textup{(a)}]
$H_\Ec$ is a dense subspace of $H.$

\smallskip
\item[\textup{(b)}]
$H_\Ec$ is a Banach space with respect to $\norm{\cdot}_{1,\Ec}.$

\smallskip
\item[\textup{(c)}]
If $x \in H_\Ec$, then
$x = \sum \, \ip{x}{e_n} \, \te_n,$
with convergence of this series in both of the norms
$\norm{\cdot}$ and $\norm{\cdot}_{1,\Ec}.$

\smallskip
\item[\textup{(d)}]
$H_\Ec$ is topologically isomorphic to $\ell^1$
via the mapping $Tx = \bigset{\ip{x}{e_n}},$ for $x \in H_\Ec.$

\smallskip
\item[\textup{(e)}]
$H_\Ec$ is a meager (first category) subset of $H.$
\qeddef
\end{enumerate}
\end{proposition}

If we replace Riesz bases by frames in Proposition \ref{iso_ell1},
then it is still true that $H_\Ec$ is a Banach space.
However, the mapping $T$ defined in statement~(d)
will not be surjective in general.
In fact, $T$ is surjective if and only if $\Ec$ is a Riesz basis.
For related results in this direction, see \cite{HLT11}.

Next we prove that the closure of an $\ell^1$-bounded set is
$\ell^1$-bounded (compare \cite[Lem.~D.2]{HH13}).

\begin{lemma} \label{ell1_closure}
Assume $M \subseteq H$ is $\ell^1$-bounded with respect to a
Riesz basis $\Ec = \set{e_n}_{n \in \N},$
and set $R = \sup_{x \in M} \norm{x}_{1,\Ec}.$
Then the closure $\overline{M}$ of $M$ is also $\ell^1$-bounded, and
$\sup_{x \in \overline{M}} \norm{x}_{1,\Ec} = R.$
\end{lemma}
\begin{proof}
Take $x \in\overline{M},$ and let $\set{x_k}$ be any sequence of elements
of $M$ that converges to $x$ in the norm of $H.$
Using Fatou's Lemma for series, we compute that
\begin{align*}
\norm{x}_{1,\Ec}
\Eq \sum_n |\ip{x}{e_n}|
& \Eq \sum_n \lim_{k \to \infty} |\ip{x_k}{e_n}|
\Le \liminf_{k \to \infty} \sum_n |\ip{x_k}{e_n}|
\Le R. \qedhere
\end{align*}
\end{proof}

\smallskip
\subsection{Minimal and Maximal $H_\Ec$ Spaces}

Does there exist a Riesz basis $\Ec$ such that
$H_\Ec \subseteq H_\Rc$ for every other Riesz basis $\Rc$?
Does there exist a Riesz basis $\Ec$ such that
$H_\Ec \supseteq H_\Rc$ for every other Riesz basis $\Rc?$
We answer both questions negatively below.

\begin{lemma} \label{inf_ell1_norm}
If $x \ne 0,$ then there is a Riesz basis $\Ec$ for $H$
such that $\norm{x}_{1,\Ec} = \infty.$
\end{lemma}
\begin{proof}
Without loss of generality, we can assume that $\norm{x} = 1.$
Let $M = \text{span}\set{x},$ and fix any sequence of scalars $(a_n)$
such that $\sum |a_n| = \infty$ and $\sum |a_n|^2 = 1.$
If we set $x_n = a_n x$ then
$\Fc = \set{x_n}$ is a frame for $M,$ and
$\norm{x}_{1,\Fc} = \infty.$

By Theorem \ref{Nai_Du}, there exists a Hilbert space $K \supseteq M$
and a Riesz basis $\Ec = \set{e_n}$ for $K$ such that $P_M e_n = x_n$
for every $n.$
Further, $K$ can be constructed to have the form $K = M \times N^{\perp},$
where $N$ is a closed subspace of $\ell^2$ topologically isomorphic to $M$
(and we identify $M$ with $M \times \set{0}$ in $K$).
We have
\begin{align*}
\sum_n |\ip{x}{e_n}|
& \Eq \sum_n |\ip{x}{P_M e_n}| 
\Eq \sum_n |\ip{x}{x_n}|
\Eq \infty.
\end{align*}

Now let $\set{y_n}$ and $\set{z_n}$ be orthonormal bases for
$M^\perp \subseteq H$ and $N^\perp \subseteq \ell^2,$ respectively.
If $y \in K,$ then
$y=\Bigparen{cx,\sum c_n z_n}$ for a unique choice of scalar $c$
and square-summable sequence $(c_n).$
Define
$$Ly
\Eq L\biggparen{cx,\sum_n c_n z_n}
\Eq cx \Plus \sum_n c_n y_n.$$
Then $L$ is a unitary mapping from $K$ onto $H,$ and $Lx = x.$
Finally, $\set{Le_n}$ is a Riesz basis for $H$ since $L$ is unitary, and
\begin{align*}
\sum_n |\ip{x}{Le_n}|
& \Eq \sum_n |\ip{L^{-1}x}{e_n}|
\Eq \sum_n |\ip{x}{e_n}|
\Eq \infty. \qedhere
\end{align*}
\end{proof}

\begin{corollary}
There does not exist a Riesz basis $\Ec$ for $H$ such that
$H_\Ec \subseteq H_\Rc$ for every Riesz basis $\Rc$ for $H.$
\qeddef\end{corollary}

\begin{remark} \label{HF_properly_contained}
A natural weaker follow-up question is:
For any Riesz basis $\Ec,$ can we find another Riesz basis $\Rc$
such that $H_\Ec \subsetneq H_\Rc$?
We can easily find a different Riesz basis $\Rc$ such that $H_\Ec = H_\Rc$
simply by rescaling $\Ec.$
However, it is unclear to the authors whether we can find a Riesz basis $\Rc$
such that $H_\Rc$ ``properly" contains $H_\Ec.$
\qeddef\end{remark}

\begin{proposition}
There does not exist a Riesz basis $\Ec$ such that
$H_\Ec \supseteq H_\Rc$ for every Riesz basis $\Rc$ for $H.$
\end{proposition}
\begin{proof}
Let $\Ec = \set{e_n}$ be an orthonormal basis for $H.$
Define $y_0 = \frac{\sqrt6}{\pi} \sum \frac{1}{n} e_n,$
so $\norm{y_0} = 1.$
Let $\set{y_n}_{n \in \N}$ be an orthonormal basis for $\set{y_0}^\perp.$
Then $\Rc = \set{y_n}_{n \ge 0}$ is an orthonormal basis for $H,$
and $\norm{y_0}_{1,\Ec} = \infty$
but $\norm{y_0}_{1,\Rc} = 1.$
Therefore $H_\Rc \nsubseteq H_\Ec.$
\end{proof}

\smallskip
\subsection{Intersections of $H_\Ec$ Spaces}
We consider how large or small the intersection of two $H_\Ec$ spaces
corresponding to Riesz bases can be.

The intersection of any two $H_\Ec$ spaces contains at least the zero element,
and we can easily find two different Riesz bases such that
$H_\Ec \cap H_\Rc = H_\Ec$ simply by letting $\Rc$ be a rescaling of $\Ec.$
Can we find two Riesz bases such that the intersection of their
corresponding $H_\Ec$ spaces is $\set{0}$?
Surprisingly, the answer is yes.
We note that the proof relies on Theorem \ref{Bes_decomposition},
which requires the underlying scalar field to be $\C$;
we do not know if a similar result holds over the real field.

\begin{proposition} \label{intersection_thm}
There exist Riesz bases $\Ec$ and $\Rc$ for $H$ such that
$H_\Rc \cap H_\Ec = \set{0}.$
\end{proposition}
\begin{proof}
Let $\Fc = \set{x_n}$ be a frame such that $H_\Fc = \set{0}$
(see Example \ref{empty_HF_space}).
By Theorem \ref{Bes_decomposition}, we can find two Riesz bases
$\Ec = \set{e_n}$ and $\Rc = \set{r_n}$ such that $x_n = e_n + r_n$
for every $n.$
Consequently, if $x \in \Ec \cap \Rc,$ then
\begin{align*}
\norm{x}_{1,\Fc}
& \Eq \sum_n |\ip{x}{x_n}|
\Le \sum_n |\ip{x}{e_n}| \Plus \sum_n |\ip{x}{r_n}|
\Eq \norm{x}_{1,\Ec} \Plus \norm{x}_{1,\Rc}
\Lt \infty.
\end{align*}
Therefore $x \in H_\Fc = \set{0},$ so we conclude that
$\Ec \cap \Rc = \set{0}.$
\end{proof}

\begin{remark}
Thus, there exist two Riesz bases such that no subset of $H$ other than
$\set{0}$ can be $\ell^1$-bounded with respect to both of these bases.
One can ask if there exist two orthonormal bases for which this is true.
It is shown in \cite{Cas98} that a frame can be written as a
linear combination of two orthonormal bases if and only if it is a Riesz basis.
Therefore the argument in the proof of Theorem \ref{intersection_thm}
does not apply if $\Ec$ and $\Rc$ are both orthonormal bases.
This leaves us with an open question about orthonormal bases.
\qeddef\end{remark}

The frame constructed in Example \ref{empty_HF_space} is not
norm-bounded below.
There are significant differences in nature between frames that are
norm-bounded below and those that are not.
It is therefore natural to ask if there is a frame $\Fc$ that is
norm-bounded below and satisfies $H_\Fc = \set{0}.$
The technique used in Example \ref{empty_HF_space} is no longer
applicable in this situation, but even so we are able to construct
such a frame.

\begin{example} \label{empty_HF_space_nbb}
Let $\Gc = \set{x_n}$ be the frame from Example \ref{empty_HF_space}
that satisfies $H_\Gc = \set{0}.$
By Theorem \ref{Bes_decomposition}, there exist Riesz bases
$\set{y_n}$ and $\set{z_n}$ such that $x_n = y_n + z_n$ for every $n.$
Let $\Fc = \set{w_n}$ be an enumeration of the union of
$\set{y_n}$ and $\set{z_n}.$
This is a (redundant) frame, and it is norm-bounded below.
If $H_\Fc \ne \set{0},$ then there is some nonzero $x \in H$ such that
$\norm{x}_{1,\Fc} < \infty.$
But then, by the triangle inequality,
$$ \norm{x}_{1,\Gc}
\Eq \sum_n |\ip{x}{x_n}|
\Le \sum_n |\ip{x}{y_n}| \Plus \sum_n |\ip{x}{z_n}|
\Eq \norm{x}_{1,\Fc}
\Lt \infty,$$
which is a contradiction.
Hence $H_\Fc = \set{0}.$
\qeddef\end{example}

\smallskip
\subsection{Unions and Sums of $H_\Ec$ Spaces} \label{union_ell1_bdd_HF}
If $\Ec$ and $\Rc$ are Riesz bases, is
$H_\Ec \cup H_\Rc \subseteq H_\Gc$ for some Riesz basis $\Gc?$ 
By the triangle inequality and the fact that $0\in H_\Ec \cap H_\Rc,$
we have $H_\Ec \cup H_\Rc \subseteq H_\Gc$
if and only if $H_\Ec + H_\Rc \subseteq H_\Gc.$
Therefore considering unions is equivalent to considering sums
of $H_\Ec$ spaces.

First we give a sufficient condition.

\begin{proposition} \label{suffi_closed_union}
Let $\Ec = \set{e_n}$ and $\Rc = \set{r_n}$ be Riesz bases for $H,$
and let $\tEc = \set{\te_n}$ and $\tRc = \set{\tr_n}$ be their dual bases,
respectively.
If there is a Riesz basis $\Gc = \set{g_n}$ such that
$\tEc$ and $\tRc$ are each $\ell^1$-bounded with respect to $\Gc,$
then $H_\Ec \cup H_\Rc \subseteq  H_\Gc.$
\end{proposition}
\begin{proof}
Let $R = \sup \norm{\te_k}_{1,\Gc}.$
If $x \in H_\Ec,$  then by applying Tonelli's Theorem for series
we see that
\begin{align*}
\sum_n |\ip{x}{g_n}|
& \Eq \sum_n \, \biggabs{\biggip{\sum_k \ip{x}{e_k} \, \te_k}{g_n}}
      \\
& \Le \sum_k \sum_n |\ip{x}{e_k}| \, |\ip{\te_k}{g_n}|
      \allowdisplaybreaks \\
& \Le R \, \sum_k |\ip{x}{e_k}|
      \\
& \Eq R \, \norm{x}_{1,\Ec}
\Lt \infty.
\end{align*}
Therefore $x \in H_\Gc.$
A similar argument shows that $H_\Rc \subseteq H_\Gc.$
\end{proof}

Somewhat surprisingly, a necessary condition for the collection of
$H_\Ec$ spaces to be closed under unions is that every countable sequence
must be contained in the $H_\Ec$ space of some Riesz basis.

\begin{theorem} \label{nece_closed_union_HF}
Assume that for any Riesz bases $\Ec$ and $\Rc,$ there exists
a Riesz basis $\Gc$ such that $H_\Ec \cup H_\Rc \subseteq H_\Gc.$
Then every countable sequence $\set{x_n}$ in $H$ is contained
in $H_\Gc$ for some Riesz basis $\Gc.$
\end{theorem}
\begin{proof}
The hypotheses of the theorem are
equivalent to assuming that if $\Ec$ and $\Rc$ are any
two Riesz bases for $H,$ then there exists a Riesz basis $\Gc$
such that $H_\Ec+ H_\Rc \subseteq H_\Gc.$
By Theorem \ref{Bes_decomposition}, this further implies that
every Bessel sequence in $H$ is contained in $H_\Gc$
for some Riesz basis $\Gc.$

Suppose that there is a countable sequence $\set{x_n}$ in $H$
that is not contained in $H_\Gc$ for any Riesz basis $\Gc.$
Fix a Riesz basis $\Ec$ for $H,$ and
choose scalars $c_n,$ all nonzero, such that $R = \sum |c_n|^2 < \infty.$
Since $H_\Ec$ is dense in $H,$ for each $n$ we can find a
vector $y_n \in H_\Ec$ such that $\norm{x_n - y_n} < |c_n|.$
Then for any $x \in H$ we have
$$\sum_n |\ip{x}{x_n-y_n}|^2
\Le \norm{x}^2 \, \sum_n \norm{x_n-y_n}^2
\Le \norm{x}^2 \, \sum_n |c_n|^2
\Eq R \, \norm{x}^2.$$
Consequently $\set{x_n-y_n}$ is a Bessel sequence in $H.$
Hence there is a Riesz basis $H_\Rc$ that contains $\set{x_n-y_n}.$
But then $\set{x_n} = \set{y_n} + \set{x_n-y_n}$ is contained in
$H_\Ec + H_\Rc,$ which itself must be contained in $H_\Gc$
for some Riesz basis $\Gc.$
However, this contradicts the definition of the sequence $\set{x_n}.$
\end{proof}

The seemingly impossible implication of Theorem $\ref{nece_closed_union_HF}$
leads us to make the following conjecture.

\begin{conjecture}
The collection of $H_\Ec$ spaces is not closed under finite unions
or finite sums.
\qeddef\end{conjecture}

\section{$\ell^1$-Bounded Sets} \label{ell1_bdd_sets_sec}

We consider $\ell^1$-bounded sets in this section.

\begin{example}
If $\Ec$ is a Riesz basis for $H,$ then it has a dual basis $\tEc$
that is also a Riesz basis for $H.$
Since $\Ec$ and $\tEc$ are biorthogonal, it follows that
$\Ec$ is $\ell^1$-bounded (with respect to the dual basis $\tEc$).
\qeddef\end{example}

\smallskip
\subsection{Characterizations of $\ell^1$-Bounded Sets}

One necessary condition for a set to be $\el$-bounded is
that it be bounded with respect to the norm of $H.$

\begin{proposition} \label{ell1_bdd_norm}
If $M$ is an $\ell^1$-bounded subset of $H,$
then $M$ is bounded in the norm of $H.$
\end{proposition}
\begin{proof}
Assume that $R = \sup_{x \in M} \norm{x}_{1,\Ec}$ is finite.
Since $\Ec$ is a Riesz basis, it is a frame,
and so has finite frame bounds $A$ and $B.$
Therefore, for each $x \in M$ we have that
$$A \, \norm{x}^2
\Le \sum_n |\ip{x}{e_n}|^2
\Le \biggparen{\sum_n |\ip{x}{e_n}|}^{\!2}
\Le R^2.$$
Hence $M$ is a bounded subset of $H.$
\end{proof}

In the finite-dimensional setting,
the $\ell^1$-norm is equivalent to the $\ell^2$-norm.
This fact allows us to give a complete characterization
of $\ell^1$-bounded sets contained in finite-dimensional subspaces.

\begin{proposition} \label{ell1_bdd_findim}
If $M \subseteq H$ is contained in a finite-dimensional
subspace of $H,$
then $M$ is $\ell^1$-bounded if and only if $M$ is bounded in the norm of $H.$
\end{proposition}
\begin{proof}
$(\Rightarrow)$ This follows from Proposition \ref{ell1_bdd_norm}.

\medskip
$(\Leftarrow)$
Assume that $M$ is contained in a subspace $S$ of dimension $d,$
and that $R = \sup_{x \in M} \norm{x}$ is finite.
Let $\set{e_1,\dots,e_d}$ be an orthonormal basis for $S,$
and let $\set{e_n}_{n=d+1}^\infty$ be an orthonormal basis for $S^\perp.$
Then $\Ec = \set{e_n}_{n \in \N}$ is an orthonormal basis for $H,$ and
\begin{align*}
\sup_{x \in M} \, \sum_n |\ip{x}{e_n}|
& \Eq \sup_{x \in M} \, \sum_{n=1}^d |\ip{x}{e_n}|
\Le d^{1/2} \, \sup_{x \in M} \biggparen{\sum_{n=1}^d |\ip{x}{e_n}|^2}^{\!1/2}
\Le d^{1/2} R,
\end{align*}
so $M$ is $\ell^1$-bounded.
\end{proof}

\smallskip
\subsection{Sets That Are Not $\ell^1$-Bounded}
First we show that Proposition \ref{ell1_bdd_findim} cannot be extended
beyond the finite-dimensional setting.

\begin{example}
For each Riesz basis $\Ec = \set{e_n}$ for $H,$
let $\tEc = \set{\te_n}$ be its dual basis, and set
$x_\Ec = \sum \, (1/n) \, \te_n.$
Then
$$\Cc
\Eq \biggset{\frac{x_\Ec}{\norm{x_\Ec}} \,:\, \Ec \text{ is a Riesz basis}}$$
is a bounded subset of $H,$
but for each Riesz basis $\Ec$ there is an element $x$ of $\Cc$
such that $\norm{x}_{1,\Ec} = \infty.$
Therefore $\Cc$ is bounded but not $\ell^1$-bounded.
\qeddef\end{example}

The set $\Cc$ constructed above is uncountable.
Next we construct a countable set that is bounded but not $\ell^1$-bounded.

\begin{example} \label{countable_non_ell1_set}
Since $H$ is separable, it contains a countable subset, say $\Dc,$
that is dense in the unit disk in $H.$
If $\Dc$ were $\ell^1$-bounded, then Lemma \ref{ell1_closure}
would imply that the unit disk in $H$ is $\ell^1$-bounded,
which is false.
Therefore $\Dc$ is bounded but not $\ell^1$-bounded.
\qeddef\end{example}

Combining Corollary \ref{countable_non_ell1_set}
and Lemma \ref{ell1_closure} gives us the following concrete example.
Here, $c_{00}$ denotes the set of sequences that contain
at most finitely many nonzero components.

\begin{corollary} \label{c00_non_ell1_bdd}
$\set{x \in c_{00} : \norm{x}_2 = 1}$ is not $\ell^1$-bounded in $\ell^2.$
\qeddef\end{corollary}

The sets that we have constructed so far that are not $\ell^1$-bounded
are dense subsets of $H.$
We will construct a subset that is not $\ell^1$-bounded but also not dense.
We will need the following lemma (see \cite[Cor.~5.9]{Car04}).

\begin{lemma} \label{singular_Lp}
Fix $1 \le p,$ $r < \infty$ with $p \ne r.$
If $T \colon \ell^p \to \ell^r$ is is a bounded linear operator,
then there do not exist closed infinite-dimensional subspaces
$S \subseteq \ell^p$ and $M \subseteq \ell^r$ such that
$T \colon S \to M$ is a topological isomorphism.
\qeddef\end{lemma}

The main ingredient in our construction
is a special case of \cite[Thm.~3.5]{HLT11}.
We include the proof for completeness.

\begin{theorem} \label{HF_nequal_H}
If $\Fc = \set{x_n}$ is a frame for $H,$ then $H_\Fc \ne H.$
\end{theorem}
\begin{proof}\,
Let $C \colon H \to \ell^2$ be the \emph{analysis operator} for $\Fc,$
defined by $Cx = \set{\ip{x}{x_n}}$ for $x \in H.$
Then $C(H_\Fc)$ is a subspace of $\ell^1,$
and we will show that is an infinite-dimensional closed subspace of $\ell^1.$

Since $C$ is injective, we know that $C(H_\Fc)$ is infinite-dimensional.
To show that it is closed, suppose that vectors
$z_n = Cy_n \in C(H_\Fc)$ and $z \in \ell^1$
are such that $z_n \to z$ in $\ell^1$-norm.
Then $Cy_n = z_n \to z$ in $\ell^2$-norm.
But $\range(C)$ is a closed subspace of $\ell^2$ since $\Fc$
is a frame (see \cite[Thm.~8.29]{Hei11}).
Therefore $z = Cy$ for some $y \in H.$
But $z \in \ell^1,$ so $y \in H_\Fc$ and hence $z = Cy \in C(H_\Fc).$
Consequently $C(H_\Fc)$ is closed in $\ell^1.$

Now consider the inclusion map $i : \ell^1 \to \ell^2.$
If $H_\Fc = H,$ then $i$ is bounded linear injective map of
the closed subspace $C(H_\Fc)$ of $\ell^1$
onto the closed subspace $i(C(H_\Fc)) = C(H) = \range(C)$ of $\ell^2.$
The Inverse Mapping Theorem therefore implies that $i$ is a
topological isomorphism that maps $C(H_\Fc)$ onto $\range(C),$
which contradicts Lemma \ref{singular_Lp}.
\end{proof}

Now we can exhibit sets that are not $\ell^1$-frame-bounded
(hence not $\ell^1$-bounded) that are not dense or even complete in $H.$

\begin{corollary} \label{ud_non_framebdd}
The unit disk $D_H = \set{x \in H : \norm{x}=1}$
is not $\ell^1$-frame-bounded in $H.$

More generally, if $S$ is any infinite-dimensional closed subspace of $H,$
then the unit disk $D_S = \set{x \in S : \norm{x}=1}$ in $S$
is not $\ell^1$-frame-bounded in $H.$
\end{corollary}
\begin{proof}
If $\Fc$ is a frame for $H$ then Proposition \ref{HF_nequal_H}
implies that there is a vector $x \in H \setminus H_\Fc.$
Consequently $y = x/\norm{x}$ belongs to the unit disk,
yet $\norm{y}_{1,\Fc} = \infty.$

Now suppose there is an infinite-dimensional closed subspace
$S$ of $H$ for which $D_S$ is $\ell^1$-frame-bounded.
Then there is a frame $\Fc = \set{x_n}$ for $H$ such that
$\sup_{x \in D_S} \sum |\ip{x}{x_n}| < \infty.$
It follows that
\begin{align*}
\sup_{x \in D_S} \sum |\ip{x}{P_S x_n}|
\Eq \sup_{x \in D_S} \sum |\ip{P_S x}{x_n}|
\Eq \sup_{x \in D_S} \sum |\ip{x}{x_n}|
\Lt \infty.
\end{align*}
However, $\Gc = \set{P_S x_n}$ is a frame for $S,$
so this implies that $D_S$ is $\ell^1$-frame-bounded in $S.$
This contradicts our work above.
\end{proof}

We close this section with an open question.

\begin{conjecture}
If a bounded subset $M$ of $H$ is uniformly separated, meaning that
$$\inf_{x \ne y \in M} \norm{x-y} \Ne 0,$$
then $M$ is $\ell^1$-bounded.
\qeddef\end{conjecture}

\smallskip
\subsection{Sums and Unions of $\ell^1$-Bounded Sets}
Any subset of an $\ell^1$-bounded set is $\ell^1$-bounded,
so the intersection of any two $\ell^1$-bounded sets is $\ell^1$-bounded
(and can even be empty).
However, the question of whether the union of two $\ell^1$-bounded sets
is $\ell^1$-bounded is much more subtle.
Interestingly, even if it is true that $\ell^1$-bounded sets are
closed under finite unions, it does not follow that $H_\Ec$ sets
are closed under unions, nor vice versa.

\begin{lemma} \label{sums_iff_unions}
The collection of $\el$-bounded sets is closed under finite unions
if and only if it is closed under finite sums.
\end{lemma}
\begin{proof} 
Let $M$ and $N$ be two $\el$-bounded subsets of $H.$

\medskip
($\Rightarrow$)
Assume that $M \cup N$ is $\el$-bounded.
Then there is a Riesz basis $\Ec = \set{e_n}$ such that
$R = \sup_{x \in M \cup N} \sum |\ip{x}{e_n}| < \infty.$
Consequently, if $x \in M$ and $y \in N$ then
$$\sum_n |\ip{x+y}{e_n}|
\Le \sum_n |\ip{x}{e_n}| \Plus \sum_n |\ip{y}{e_n}|
\Le 2R,$$
so $M+N$ is $\el$-bounded.

\medskip
($\Leftarrow$)
If $M$ and $N$ are $\el$-bounded, then so are
$M \cup\set{0}$ and $N \cup\set{0}.$
Consequently, if $\ell^1$-bounded sets are closed under finite sums then
$M \cup N \subseteq (M \cup \set{0}) + (N \cup\set{0})$ is $\el$-bounded.
\end{proof}

Now we show that several unexpected results follow
if the collection of $\el$-bounded sets is closed under finite unions.
In the statement of part~(d) of the next theorem, a sequence
$\set{x_n}$ is said to be \emph{unconditionally summable} if
the series $\sum x_n$ converges unconditionally in $H$
(that is, it converges regardless of the ordering of the index set).

\begin{theorem} \label{closed_union_implications}
If the collection of $\el$-bounded subsets of $H$
is closed under finite unions, then the following statements hold.

\begin{enumerate}
\item[\textup{(a)}]
Every Bessel sequence in $H$ (and hence every frame) is $\ell^1$-bounded.

\medskip
\item[\textup{(b)}]
$\ell^1$-frame-boundedness is equivalent to $\ell^1$-boundedness in $H.$

\medskip
\item[\textup{(c)}]
If $\Ec$ is a Riesz basis for $H,$ there then exists another
Riesz basis $\Rc$ such that $H_\Ec \subsetneq H_\Rc.$

\medskip
\item[\textup{(d)}]
Every unconditionally summable sequence is $\el$-bounded.
\end{enumerate}
\end{theorem}
\begin{proof}
Assume that the collection of $\el$-bounded sets
is closed under finite unions.

\medskip
(a) A Riesz basis is $\ell^1$-bounded with respect to its biorthogonal system.
By Theorem \ref{Bes_decomposition}, every Bessel sequence is contained
in the sum of two Riesz bases.
Hence every Bessel sequence is a subset of an $\el$-bounded set,
and therefore is itself $\el$-bounded.

\medskip
(b) Assume that $M$ is $\ell^1$-bounded set with respect to some
Riesz basis $\Ec.$
Since every Riesz basis is a frame, it follows that $M$
is $\ell^1$-frame-bounded.

Next, assume that $M$ is $\ell^1$-frame-bounded with respect to
some frame $\Fc = \set{x_n}.$
Then there exists a frame $\set{y_n}$ (such as the canonical dual frame)
such that $x = \sum \, \ip{x}{x_k} \, y_k$ for every $x \in H.$
Now, statement~(a) implies that $\set{y_n}$ is $\ell^1$-bounded
with respect to some Riesz basis $\Ec = \set{e_n},$
so $R = \sup_k \norm{y_k}_{1,\Ec} < \infty.$
Then for each $x \in M$ we use Tonelli's Theorem to compute that
\begin{align*}
\sup_{x \in M} \norm{x}_{1,\Ec}
\Eq \sup_{x \in M} \biggparen{\sum_n |\ip{x}{e_n}|}
& \Eq \sup_{x \in M} \biggparen{\sum_n \,
      \biggabs{\biggip{\sum_k \ip{x}{x_k} \, y_k}{e_n}}}
      \\
& \Le \sup_{x \in M}
      \biggparen{\sum_n \sum_k |\ip{x}{x_k}| \, |\ip{y_k}{e_n}|}
      \allowdisplaybreaks \\
& \Eq \sup_{x \in M} \biggparen{\sum_k |\ip{x}{x_k}| \, \norm{y_k}_{1,\Ec}}
      \\[1 \jot]
& \Le R \, \sup_{x \in M} \norm{x}_{1,\Fc}
\Lt \infty.
\end{align*}
Thus $M$ is $\ell^1$-bounded with respect to the Riesz basis $\Ec.$

\medskip
(c) Suppose to the contrary that there exists a Riesz basis
$\Ec = \set{e_n}_{n \in \N}$ such that if $H_\Ec \subseteq H_\Rc$
for some Riesz basis $\Rc,$ then $H_\Ec = H_\Rc.$ 

Let $\tEc = \set{\te_n}_{n \in \N}$ be the dual basis to $\Ec,$
and let $x_0 = \sum \, (1/n) \, \te_n,$
which does not belong to $H_\Ec.$
Each of $\set{\te_n}_{n \in \N}$ and $\set{x_0,0}$ are $\el$-bounded sets,
and we let
$$S \Eq \set{\te_n}_{n \in \N} \Plus \set{x_0,0}
\Eq \set{\te_n}_{n \in \N} \,\cup\, \set{x_0 + \te_n}_{n \in \N}.$$
By hypothesis, $S$ is $\ell^1$-bounded with respect to some
Riesz basis $\Rc = \set{r_n}.$
Let $R = \sup_{x \in S} \norm{x}_{1,\Rc}.$

Arguing similarly to the proof of statement (b), we see
that for every $x \in H_\Ec$ we have
\begin{align*}
\norm{x}_{1,\Rc}
\Eq \sum_n |\ip{x}{r_n}|
& \Eq \sum_n \, \biggabs{\biggip{\sum_k \ip{x}{e_k} \, \te_k}{r_n}}
      \\
& \Le \sum_n \sum_k \, |\ip{x}{e_k}| \, |\ip{\te_k}{r_n}|
      \allowdisplaybreaks \\
& \Eq \sum_k \, |\ip{x}{e_k}| \, \norm{\te_k}_{1,\Rc}
     \\[1 \jot]
& \Le R \, \norm{x}_{1,\Ec}
\Lt \infty.
\end{align*}
Thus $x \in H_\Rc,$ so $H_\Ec \subseteq H_\Rc$
and hence $H_\Ec = H_\Rc.$
However, $\te_1 + x_0 \in H_\Rc \setminus H_\Ec,$ which is a contradiction.

\medskip
(d) If $\set{x_n}$ is an unconditionally summable sequence in $H,$
then $\sum \norm{x_n}^2 < \infty$
by Orlicz's Theorem (see \cite[Thm.3.16]{Hei11}).
Hence
$\sum |\ip{x}{x_n}|^2 \Le \norm{x}^2 \sum \norm{x_n}^2 < \infty$
for every $x \in H.$
The Uniform Boundedness Principle therefore implies that $\set{x_n}$
is a Bessel sequence, so the result follows from statement~(a).
\end{proof}

For some related results on unconditionally summable sequences,
we refer to \cite{FGP17}.

Given such strong implications, we make the following conjecture.

\begin{conjecture}
The collection of $\ell^1$-bounded sets is not closed
under finite unions or sums.
\end{conjecture}

\smallskip
\subsection{The Relation Between $\ell^1$-Boundedness and
$\ell^1$-Frame Boundedness}
We saw in Theorem \ref{closed_union_implications} that
$\ell^1$-frame-boundedness is equivalent to $\ell^1$-boundedness
if the class of $\ell^1$-bounded sets is closed under finite sums.
Clearly, $\ell^1$-boundedness implies $\ell^1$-frame-boundedness
since Riesz bases are non-overcomplete frames.

In the other direction, one might expect that every frame contains
a Riesz basis because of its overcompleteness, and conclude that
$\ell^1$-frame-boundedness is indeed equivalent to $\ell^1$-boundedness.
However, this is not the case.
It is easy to construct a frame that is not norm-bounded below that does not
contain Riesz basis (for example, see \cite[Example 8.6]{Hei11}).
Casazza and Christensen \cite{CC98} constructed a frame that is
norm-bounded below but does not contain any Riesz bases
(or even any Schauder bases, see \cite{CC98b}).
In this section, we study the relation between $\ell^1$-boundedness
and $\ell^1$-frame-boundedness in more detail.

First, we give another characterization of when
$\ell^1$-boundedness and $\ell^1$-frame-boundedness are equivalent.

\begin{theorem} \label{equi_frame_Riesz}
If $M \subseteq H,$ then following statements are equivalent.

\begin{enumerate}
\item[\textup{(a)}]
If $H$ embeds into a separable Hilbert space $K$
and $M$ is $\ell^1$-bounded in $K,$
then $M$ is $\ell^1$-bounded in $H.$

\medskip
\item[\textup{(b)}]
$M$ is $\ell^1$-frame-bounded in $H$ if and only if
$M$ is $\ell^1$-bounded in $H.$
\end{enumerate}
\end{theorem}
\begin{proof}
(a) $\Rightarrow$ (b).
Assume that statement~(a) holds.
If $M$ is $\ell^1$-bounded then certainly $M$ is $\ell^1$-frame-bounded,
so assume that $M$ is $\ell^1$-frame-bounded with respect to a frame
$\Fc = \set{x_n}$ for $H.$
By Theorem \ref{Nai_Du}, there is a Hilbert space $K$
and a Riesz basis $\Ec = \set{e_n}$ for $K$ such that
$H$ embeds into $K$ and $P_H e_n = x_n$ for every $n.$
Therefore $M$ is $\ell^1$-bounded in $K$ since
$\sum |\ip{x}{e_n}| \Eq \sum |\ip{x}{P_H e_n}| < \infty.$
Statement~(a) therefore implies that $M$ is $\ell^1$-bounded in $H.$

\medskip
(b) $\Rightarrow$ (a).
Suppose that statement~(b) holds.
Assume that $H$ embeds into a separable Hilbert space $K,$
and that $M$ is $\ell^1$-bounded with respect to some
Riesz basis $\Ec = \set{e_n}$ for $K.$
If we set $x_n = P_H e_n,$ then $\Fc = \set{x_n}$ is a frame for $H.$
For $x \in M$ we have
$$\norm{x}_{1,\Fc}
\Eq \sum_n |\ip{x}{x_n}|
\Eq \sum_n |\ip{x}{P_H e_n}|
\Eq \sum_n |\ip{P_H x}{e_n}|
\Eq \sum_n |\ip{x}{e_n}|
\Eq \norm{x}_{1,\Ec}.$$
Therefore $M$ is $\ell^1$-frame-bounded in $H.$
Consequently statement~(b) implies that $M$ is $\ell^1$-bounded in $H.$
\end{proof}

\smallskip
\subsection{Embeddings and $\ell^1$-Boundedness}
Theorem \ref{equi_frame_Riesz} suggests that we should consider how
$\ell^1$-boundedness relates to embeddings.
One direction is clear.

\begin{proposition} \label{embedding_prop}
Assume that $H$ embeds into a separable Hilbert space $K.$
If $M \subseteq H$ is $\ell^1$-bounded in $H,$
then $M$ is $\ell^1$-bounded in $K.$
\end{proposition}
\begin{proof}
Assume $M$ is $\ell^1$-bounded with respect to some Riesz basis $\Ec$
for $H,$ and let $\Rc$ be a Riesz basis for $H^\perp$ in $K.$
Then $\Ec \cup \Rc$ is a Riesz basis for $K,$
and $M$ is $\ell^1$-bounded with respect to $\Ec \cup \Rc.$
\end{proof}

We do not know whether the converse of Proposition \ref{embedding_prop}
holds in general.
We give proofs for two special cases.
First, we assume that the codimension of $\clspan(M)$ in $H$ is infinite.

\begin{proposition} \label{l1_bd_lar_imply_l1_bd_sma}
Assume that $H$ embeds into a separable Hilbert space $K,$
and $M \subseteq H$ is $\ell^1$-bounded in $K.$
Let $M^\perp_H$ be the orthogonal complement of $M$ in $H.$
If $\text{dim}(M^\perp_H) = \infty,$
then $M$ is $\ell^1$-bounded in $H.$
\end{proposition}
\begin{proof}
Let $S = \clspan{(M)}.$
Let $S_H^\perp$ denote the orthogonal complement of $S$ in $H,$
and let $S_K^\perp$ denote the orthogonal complement of $S$ in $K.$
By hypothesis, $S_H^\perp$ and $S_K^\perp$ are each separable,
infinite-dimensional Hilbert spaces.
Let $U \colon S_K^\perp \to S_H^\perp$ be a unitary mapping.
Given $x \in K,$ write
$x
= x_S \Plus x_{S_K^\perp},$ 
where $x_S \in S$ and $x_{S_K^\perp} \in S_K^\perp,$
and define $T \colon K \to H$ by
$$Tx
\Eq x_S \Plus U(x_{S^\perp_K}),
\qquad\text{for } x \in K.$$
Then $T$ is a unitary map that fixes $S$ (and hence $M$).
Suppose that $M$ is $\ell^1$-bounded with respect to some
Riesz basis $\Ec = \set{e_n}$ for $K.$
Then  $T(\Ec) = \set{Te_n}$ is a Riesz basis for $H,$ and we have
\begin{align*}
\sup_{x \in M} \norm{x}_{1,T(\Ec)}
\Eq \sup_{x \in M} \sum_n |\ip{x}{Te_n}|
& \Eq \sup_{x \in M} \sum_n |\ip{Tx}{Te_n}|
\Eq \sup_{x \in M} \sum_n |\ip{x}{e_n}|
\Lt \infty.
\end{align*}
Therefore $M$ is $\ell^1$-bounded in $H.$
\end{proof}

\begin{corollary}
Let $M$ be a subset of $H$ such that
$\text{dim}(M^\perp) = \infty.$
Then $M$ is $\ell^1$-bounded in $H$
if and only if it is $\ell^1$-frame-bounded in $H.$
\end{corollary}
\begin{proof}
This follows by applying
Proposition \ref{l1_bd_lar_imply_l1_bd_sma}
and Theorem \ref{equi_frame_Riesz}.
\end{proof}

Our next goal is to show that $\ell^1$-boundedness
in a larger Hilbert space $K$ implies $\ell^1$-boundedness in $H$
when the codimension of $H$ in $K$ is finite.
To prove this we will need the following lemmas.

\begin{lemma} \label{criterion_Reisz_basis}
Let $\Ec = \set{e_n}$ be a Riesz basis for $H,$
with frame bounds $A$ and $B.$
If $\set{x_n}$ satisfies
$\sum \norm{e_n-x_n}^2 < A,$
then $\set{x_n}$ is a Riesz basis for $H.$
\end{lemma}
\begin{proof}
The hypotheses imply that $\set{e_n - x_n}$ is a Bessel sequence in $H$
with upper frame bound strictly less than $A.$
The result therefore follows from \cite[Cor.~22.1.5]{Chr16}.
\end{proof}

\begin{lemma} \label{Ortho_frame_twonorm}
Let $\Fc = \set{x_n}$ be a frame for $H,$ with frame bounds $A$ and $B.$
If $S$ is a closed subspace of $H,$ then
$\sum \norm{P_S x_n}^2 \Le B \, \text{dim}(S).$
\end{lemma}
\begin{proof}
Let $\set{e_k}_{k \in J}$ be an orthonormal basis for $S,$
and then extend this to an orthonormal basis $\set{e_k}_{k \in \N}$ for $H.$
Then $P_S e_k = e_k$ for $k \in J$ and $P_S e_k = 0$ for $k \notin J,$ so
\begin{align*}
\sum_n \norm{P_S x_n}^2
& \Eq \sum_n \sum_k |\ip{P_Sx_n}{e_k}|^2 \\
& \Eq \sum_{k \in J} \sum_n |\ip{x_n}{P_S e_k}|^2
\Le B \sum_{k \in J} \norm{P_S e_k}^2
\Eq B \, \text{dim}(S). \qedhere
\end{align*}
\end{proof}

The proof of our next theorem was inspired by the results of
Holub in \cite{Hol94}.

\begin{theorem}
\label{Ortho_compl_H_findim}
Fix $M \subseteq H,$ and assume that $H$ embeds into a
separable Hilbert space $K.$
Let $H^\perp$ denote the orthogonal complement of $H$ in $K.$
If $M$ is $\ell^1$-bounded in $K$ and $\text{dim}(H^\perp) < \infty,$
then $M$ is $\ell^1$-bounded in $H.$
\end{theorem}
\begin{proof}
Assume that $M$ is $\ell^1$-bounded in $K$ with respect to a
Riesz basis $\Ec = \set{e_n},$
and let $A$ and $B$ be frame bounds for $\Ec.$
We have $\sum \norm{P_{H^\perp} e_n}^2 < \infty$
by Lemma \ref{Ortho_frame_twonorm}.
Choose $n_0$ large enough such that
$\sum_{n > n_0} \norm{P_{H^\perp}e_n}^2 < A$
and define
$z_n = e_n$ for $1 \le n \le n_0$ and
$z_n = P_He_n$ for $n > n_0.$
Then
$$\sum_n \norm{e_n-z_n}^2
\Eq \sum_{n > n_0} \norm{P_{H^\perp}e_n}^2
\Lt A,$$
so $\set{z_n}$ is a Riesz basis for $K$ by Lemma \ref{criterion_Reisz_basis}.
Therefore, $\set{z_n}_{n > n_0} = \set{P_H e_n}_{n > n_0}$
is a \emph{Riesz sequence} in $H$;
that is, it is a Riesz basis for its closed span in $H.$
Moreover, it only requires finitely many elements to extend $\set{z_n}$
to a Riesz basis for $H.$
Let $y_1,\dots,y_k$ be elements of $H$ such that
$\set{y_n}_{1 \le n \le k} \cup \set{z_n}_{n > n_0}$
is a Riesz basis for $H.$
Then
\begin{align*}
\sup_{x \in M} \biggparen{\sum_{n=1}^k |\ip{x}{y_n}| \Plus
      \sum_{n > n_0} |\ip{x}{z_n}|}
& \Eq \sup_{x \in M} \biggparen{\sum_{n=1}^k |\ip{x}{y_n}| \Plus
      \sum_{n > n_0} |\ip{P_Hx}{e_n}|}
      \\
& \Eq \sup_{x \in M} \biggparen{\sum_{n=1}^k |\ip{x}{y_n}| \Plus
      \sum_{n > n_0} |\ip{x}{e_n}|}
      \\
& \Le \sup_{x \in M} \biggparen{\sum_{n=1}^k \norm{x} \, \norm{y_n} \Plus
      \sum_{n=1}^\infty |\ip{x}{e_n}|}
\Lt \infty,
\end{align*}
where at the last step we have used the facts that $M$ is bounded in norm
and that $M$ is $\ell^1$-bounded in $K.$
We conclude that $M$ is $\ell^1$-bounded in $H.$
\end{proof}

\section{$p$-Convergent Frames} \label{absolutelyconvergent_sec}

Let $\set{x_n}$ be a frame for $H$ and let $\set{y_n}$
be an associated alternative dual.
The study of convergence of the frame expansions
\begin{equation} \label{frameexpansion_eq}
x \Eq \sum_n \, \ip{x}{y_n} \, x_n
\end{equation}
has a long history.
If $\set{y_n}$ is a frame, then the series
$\sum \, \ip{x}{y_n} \, x_n$ converges unconditionally for every $x \in H.$
However, there exist frames $\set{x_n}$ and alternative duals $\set{y_n}$
that are not frames for which the frame expansion
in equation \eqref{frameexpansion_eq} converges conditionally
for some $x \in H$ (see \cite{HY22}, \cite{ST11}).
Excluding the case where $\set{x_n}$ contains infinitely zeros,
it was shown in \cite{HY22} that if the frame expansion
in equation \eqref{frameexpansion_eq}
converges unconditionally for every $x \in H$ for all
associated alternative duals $\set{y_n},$
then $\set{x_n}$ is a Riesz basis plus finitely many elements
(also known as a \emph{near-Riesz basis}).
In this section, instead of unconditional convergence,
we will consider absolute convergence of frame expansions.
Specifically, we ask whether there exists a frame and an
associated alternative dual such that the frame expansions
converge absolutely for every $x \in H.$
Using the language of $\ell^1$-boundedness, this is equivalent to asking
if there exists a frame $\set{x_n}$ and
associated alternative dual $\set{y_n}$ such that
$H_\Yc = H$ when $\Yc$ is the renormalized system
$\Yc = \set{\norm{x_n} \, y_n}.$

\begin{definition}[$p$-Convergent Frame]
Let $\set{x_n}$ be a frame for $H.$
If $1 \le p \le 2$ then we say $\set{x_n}$ is a \emph{$p$-convergent frame}
if there exists an alternative dual $\set{y_n}$ for $\set{x_n}$ such that
$$\sum_n \norm{\ip{x}{y_n} \, x_n}^p \Lt \infty,
\quad\text{for every } x \in H.$$
We say that a $1$-convergent frame is an
\emph{absolutely convergent frame}.
\qeddef\end{definition}

If $\set{x_n}$ is a frame, then it has at least one
alternative dual frame $\set{y_n}$ (the canonical dual frame).
In this case the frame expansions in equation \eqref{frameexpansion_eq}
converge unconditionally for every $x,$ because $\set{y_n}$ is a frame.
Orlicz's Theorem therefore implies that
$\sum \norm{\ip{x}{y_n} \, x_n}^2 < \infty$ for every $x$
(see \cite[Thm.~3.16]{Hei11}).
Hence every frame is a $2$-convergent frame.

We show next that if $1 < p < 2$ then there exists a $p$-convergent frame.

\begin{example} \label{exist_p_abs_frames}
Let $\set{e_n}$ be an orthonormal basis for $H.$
Fix $1<p<2,$ and let $k\in \N$ be large enough that
$k > \frac{2-p}{4(p-1)}.$
Let
$$\set{x_n}
\Eq \set{e_1} \,\cup\, \set{2^{-k}e_2}_{j=1}^{2^{2k}} \,\cup\,
    \set{3^{-k}e_3}_{j=1}^{3^{2k}} \,\cup\, \cdots \,\cup\,
    \set{n^{-k}e_k}_{j=1}^{n^{2k}} \,\cup\, \cdots.$$
Then $\set{x_n}$ is a Parseval frame.
Therefore its canonical dual frame is simply $\set{x_n},$
so the frame expansion of $x \in H$ is
$$x
\Eq \sum_{n=1}^\infty \ip{x}{x_n} \, x_n
\Eq \ip{x}{e_1} \Plus \sum_{j=1}^{2^{2k}} \ip{x}{2^{-k}e_2} \, 2^{-k}e_2
    \Plus \sum_{j=1}^{3^{2k}} \ip{x}{3^{-k}e_3} \, 3^{-k}e_3 \Plus \cdots.$$
Consequently, by applying Holder's inequality we see that
\begin{align*}
\sum_{n=1}^\infty \norm{\ip{x}{x_n} \, x_n}^p
& \Eq \sum_{n=1}^\infty n^{(2-2p)k} \, |\ip{x}{e_n}|^p
      \\[1 \jot]
& \Le \biggparen{\sum_{n=1}^\infty |\ip{x}{e_n}|^2}^{\!p/2} \,
      \biggparen{\sum_{n=1}^\infty
      n^{k(2-2p) \times \frac{2}{2-p}}}^{\!(2-p)/2}
      \\[1 \jot]
& \Eq \norm{x}^p \,
      \biggparen{\sum_{n=1}^\infty n^{-k \frac{4(p-1)}{2-p}}}^{\!(2-p)/2}
\Lt \infty.
\end{align*}
Therefore $\set{x_n}$ is a $p$-convergent frame.
\qeddef\end{example}

In contrast, since $H$ is not isomorphic to $\ell^p$ when $1 \le p < 2,$
there cannot exist a $p$-convergent Riesz basis for $H.$
Although the preceding example shows that $p$-convergent frames exist
for $1 < p < 2,$ we show next that there are no
absolutely convergent frames.

\begin{theorem} \label{nonexis_abs_conv_frame}
There does not exist an absolutely convergent frame for $H.$
\end{theorem}
\begin{proof}
Assume $\set{x_n}$ is an absolutely convergent frame for $H.$
Then there must exist an alternative dual $\set{y_n}$ such that
$\sum \norm{\ip{x}{y_n} \, x_n} < \infty$ for every $x,$
and consequently $H_\Yc = H$ where $\Yc = \set{\norm{x_n} \, y_n}.$
The the linear operator $C \colon H \to \ell^1$ defined by
$Cx = \bigset{\ip{x}{\norm{x_n} \, y_n}}$ satisfies
\begin{equation} \label{Cbounded_eq}
\norm{x}
\Eq \biggnorm{\sum_n \ip{x}{y_n} \, x_n}
\Le \sum_n \bigabs{\ip{x}{\norm{x_n} \, y_n}}
\Eq \norm{Cx}_1
\Lt \infty.
\end{equation}
The Banach--Steinhaus Theorem therefore implies that $C$ is bounded.
Moreover, the lower bound in equation \eqref{Cbounded_eq} implies
that $C$ is injective, and that $C(H)$ is a closed subspace of $\ell^1.$
Therefore $C$ is a topological isomorphism that maps $H$ onto $C(H).$
But $H$ is topologically isomorphic to $\ell^2,$
so this implies that $\ell^2$ is topologically isomorphic to the
closed subspace $C(H)$ of $\ell^1.$
This contradicts Lemma \ref{singular_Lp}.
\end{proof}

We also have a symmetric version of Theorem \ref{nonexis_abs_conv_frame}.

\begin{theorem} \label{dual_nonexis_abs_conv_frame}
There does not exist a frame $\set{x_n}$ for $H$ that has an
alternative dual $\set{y_n}$ such that
$$\sum_n \norm{\ip{x}{x_n} \, y_n} \Lt \infty,
\qquad\text{for every } x \in H.$$
\end{theorem}
\begin{proof}
Since absolute convergence implies unconditional convergence,
Lemma $\ref{equiv_conv_alt_syn_dual}$ implies that
$\sum \ip{x}{x_n} \, y_n = x$ for every $x \in H.$
The argument used in Theorem \ref{nonexis_abs_conv_frame}
then completes the proof.
\end{proof}

\begin{corollary} \label{non_exist_abs_frame_coro}
If $\set{x_n}$ is a frame for $H,$ then $\sum \norm{x_n}^2 = \infty.$
\end{corollary}
\begin{proof}
Suppose that $\set{x_n}$ is a frame such that
$\sum \norm{x_n}^2 < \infty,$ and let $S$ be its associated frame operator.
By assumption,
$$\sum_n \norm{S^{-1/2}x_n}^2
\Le \sum_n \norm{S^{-1/2}}^2 \, \norm{x_n}^2
\Lt \infty.$$
Since $\set{S^{-1/2}x_n}$ is a Parseval frame,
its canonical dual frame is itself.
However, if $x \in H$ then
$$\sum_n \norm{\ip{x}{S^{-1/2}x_n} \, S^{-1/2}x_n}
\Le \norm{x} \sum_n \norm{S^{-1/2}x_n}^2
\Lt \infty,$$
which implies that $\set{S^{-1/2}x_n}$ is an absolutely convergent frame.
This contradicts Theorem \ref{nonexis_abs_conv_frame}.
\end{proof}

The exponent $2$ in Corollary $\ref{non_exist_abs_frame_coro}$ is tight,
in the sense that if $p > 1$ then there exists a frame $\set{x_n}$
such that $\sum \norm{x_n}^{2p} < \infty.$
In particular, the frame constructed in Example $\ref{exist_p_abs_frames}$
satisfies $\sum \norm{x_n}^{2p} = \sum n^{-2k(p-1)},$
which is finite if choose $k$ to be larger than $\frac{1}{2(p-1)}.$

\end{document}